\documentclass[12pt,reqno]{amsart}
\usepackage{amsmath,amssymb,amsfonts,amscd,latexsym,amsthm,mathrsfs,verbatim}
\usepackage[usenames]{color}
\usepackage[unicode]{hyperref}
\usepackage{hypbmsec}
\theoremstyle{remark}

\newtheorem{example}{\bf Example}

\definecolor{orange}{rgb}{0.7,0.3,0}

\renewcommand{\d}{{\mathrm d}}
\begin{document}

\title{C\protect\lowercase{rouching} AGM, H\protect\lowercase{idden} M\protect\lowercase{odularity}}

\author{Shaun Cooper}
\address{Institute of Natural and Mathematical Sciences, Massey University\,---\,Albany, Private Bag 102904, North Shore Mail Centre, Auckland 0745, NEW ZEALAND}
\email{s.cooper@massey.ac.nz}

\author{Jes\'us Guillera}
\address{Department of Mathematics, University of Zaragoza, 50009 Zaragoza, SPAIN}
\email{jguillera@gmail.com}

\author{Armin Straub}
\address{Department of Mathematics and Statistics, University of South Alabama, 411 University Blvd N, ILB 325, Mobile, AL 36688, UNITED STATES}
\email{straub@southalabama.edu}

\author{Wadim Zudilin}
\address{School of Mathematical and Physical Sciences, The University of Newcastle, Callaghan NSW 2308, AUSTRALIA}
\email{wadim.zudilin@newcastle.edu.au}

\date{5 April 2016}

\begin{abstract}
Special arithmetic series $f(x)=\sum_{n=0}^\infty c_nx^n$, whose coefficients $c_n$ are normally given as certain binomial sums,
satisfy `self-replicating' functional identities. For example, the equation
$$
\frac1{(1+4z)^2}\,f\biggl(\frac z{(1+4z)^3}\biggr)
=\frac1{(1+2z)^2}\,f\biggl(\frac{z^2}{(1+2z)^3}\biggr)
$$
generates a modular form $f(x)$ of weight 2 and level 7, when a related modular parametrization $x=x(\tau)$ is properly chosen.
In this note we investigate the potential
of describing modular forms by such self-replicating equations as well as applications
of the equations that do not make use of the modularity. In particular, we outline a new recipe of generating AGM-type algorithms
for computing~$\pi$ and other related constants. Finally, we indicate some possibilities to extend the functional equations
to a two-variable setting.
\end{abstract}

\dedicatory{To Mourad Ismail, with warm wishes from three continents}


\thanks{The work of W.Z.\ is supported by the Australian Research Council.}

\maketitle

\section{Introduction}
\label{sec1}

Modular forms and functions form a unique enterprise in the world of special functions. They are
innocent looking $q$-series, at the same time highly structured and
possess numerous links to other parts of mathematics.
In our exposition below we will try to hide the modularity of objects under consideration as much as possible,
just pointing out related references to the literature where the modular origin is discussed in detail.

One way of thinking of modular forms is through Picard--Fuchs linear differential equations \cite[Section 5.4]{Zag08}.
The concept can be illustrated on the example
\begin{equation}
\label{eq:f7}
f_7(x)=\sum_{n=0}^\infty u_nx^n,
\end{equation}
with coefficients
\begin{equation}
\label{eq:s7}
u_n=\sum_{k=0}^n{\binom nk}^2\binom{n+k}n\binom{2k}n
=\sum_{k=0}^n(-1)^{n-k}\binom{3n+1}{n-k}{\binom{n+k}{n}}^3.
\end{equation}
The function $f_7(x)$ represents a level~7 and weight~2 modular form e.g., \cite{Coo12,Zud14}
(see also Section~\ref{sec3} below),
and satisfies a linear differential equation, which we record in the form of a recurrence relation for the coefficients
\begin{equation}
\label{u-rec}
(n+1)^3u_{n+1}=(2n+1)(13n^2+13n+4)u_n+3n(3n-1)(3n+1)u_{n-1}
\end{equation}
where $n=0,1,2,\dots$\,. The single initial condition $u_0=1$ is used to start the sequence.
This is sequence A183204 in the \emph{OEIS} \cite{Slo16}.
Recursions such as \eqref{u-rec} are famously known as Ap\'ery-like recurrence equations \cite{AvSZ11,Coo12,Zag09}.
The series defining $f_7(x)$ converges in the disc $|x|<1/27$ as $x=1/27$ is a singular point of the underlying linear differential equation
produced by the recursion.

In this note we will investigate another point of view on the modularity, using functional equations of a different nature.
The function $f(x)=f_7(x)$ defined above will be our principal example, while the special functional equation it
satisfies will be~\cite[Theorem~4.4]{CY16}
\begin{equation}
\frac1{(1+4z)^2}\,f\biggl(\frac z{(1+4z)^3}\biggr)
=\frac1{(1+2z)^2}\,f\biggl(\frac{z^2}{(1+2z)^3}\biggr),
\label{alg}
\end{equation}
valid in a suitable neighbourhood of $z=0$. Note that this equation,
together with the initial condition $f(0)=1$, uniquely determines the coefficients
in the Taylor series expansion $f(x)=\sum_{n=0}^\infty u_nx^n$.
In Section~\ref{sec2} we will address this fact and discuss some variations of Equation \eqref{alg}.
In particular, we observe that, in arithmetically interesting instances\,---\,those with connections to modular forms\,---\,the corresponding
sequences appear to satisfy (and possibly are characterized by) strong divisibility properties.
Without pretending to fully cover the topic of such functional equations,
we list several instances for other modular forms in Section~\ref{sec3} together with some details of the underlying
true modular forms and functions. An important feature of this way of thinking of modular forms
is applications. One such application\,---\,to establishing Ramanujan-type formulas for $1/\pi$\,---\,has already received attention
in the literature, e.g., see \cite{CZ15,GZ13,Zud13,Zud14}. In Section~\ref{sec4} we give some related examples for $f_7(x)$ and equation~\eqref{alg},
though concentrating mostly on another application of the special functional equations\,---\,to AGM-type algorithms.
Such algorithms for computing $\pi$ were first discovered independently by
R.~Brent \cite{Bre76} and E.~Salamin \cite{Sal76} in~1976 (see also \cite{BB16}).
A heavy use of the modularity made later, faster versions by other authors somehow sophisticated.
Our exposition in Section~\ref{sec4} follows the ideas in~\cite{Gui16} and
streamlines the production of AGM-type iterations by freeing them from the modularity argument.

Three-term recurrence relations, which are more general than the one given in~\eqref{u-rec}, naturally appear in
connection with orthogonal polynomials \cite{Ism05}\,---\,one of the topics mastered by Mourad Ismail.
The polynomials have long-existing ties with arithmetic, and the recent work \cite{CWZ13,RS13,WZ12,Zud14}
indicates some remarkable links between generating functions of classical Legendre polynomials and arithmetic hypergeometric functions.
Interestingly enough, these links allow one to produce a two-variable functional equation that may be thought of as a generalization
of~\eqref{alg}. This is the subject of Section~\ref{sec5}.

\section{Self-replication vs modularity}
\label{sec2}

The first natural question to ask is how sensitive the dependence of $f(x)$ on the parameters of equation~\eqref{alg} is.
In other words, what can be said about a solution $f(x)=f(\lambda,\mu;x)=\sum_{n=0}^\infty c_nx^n$, with $c_0=1$, to the equation
\begin{equation}
\frac1{(1+\mu z)^2}\,f\biggl(\frac z{(1+\mu z)^3}\biggr)
=\frac1{(1+\lambda z)^2}\,f\biggl(\frac{z^2}{(1+\lambda z)^3}\biggr)
\label{alg0}
\end{equation}
when $\lambda$ and $\mu$ are integers? Writing the equation as
$$
\sum_{n=0}^\infty\frac{c_nz^n}{(1+\mu z)^{3n+2}}
=\sum_{n=0}^\infty\frac{c_nz^{2n}}{(1+\lambda z)^{3n+2}}
$$
and applying the binomial theorem to expand both sides into a power series in $z$, we obtain
the recursion
\begin{equation}
c_n=\sum_{k=0}^{\lfloor n/2\rfloor}\binom{n+k+1}{3k+1}(-\lambda)^{n-2k}c_k-\sum_{k=0}^{n-1}\binom{n+2k+1}{3k+1}(-\mu)^{n-k}c_k,
\label{rec}
\end{equation}
for $n=1,2,\ldots$, which together with $c_0=1$ determines the sequence $c_n=c_n(\lambda,\mu)$.
Clearly, $c_n(\lambda,\mu)\in\mathbb Z$ whenever $\lambda,\mu\in\mathbb Z$.
Also, the recursion \eqref{rec} implies
that the functional equation \eqref{alg0} completely recovers the analytic solution $f(x)$. The latter fact allows us to think
of \eqref{alg} and, more generally, of \eqref{alg0} as \emph{self-replicating} equations that define the series and function $f(x)$.
In spite of the property $c_n\in\mathbb Z$ it remains unclear whether these sequences have some arithmetic significance in general.
In the sequel, we will refer to instances as \emph{arithmetic} if the generating function $f(x)$ is of modular origin.

The fact that the sequences arising from self-replicating equations such as \eqref{alg} are integral is useful since
this integrality is rather inaccessible from associated differential equations or recurrences
such as \eqref{u-rec}. (It remains an open problem to classify all integer
solutions to recurrences of the shape \eqref{u-rec}, though it is conjectured
that the known Ap\'ery-like sequences form a complete list).  For Ap\'ery-like sequences, the integrality usually follows
from binomial sum representations such as \eqref{eq:s7}.
Integrality can be also be deduced by exhibiting a modular parametrization of the differential equation;
however, we do not know \emph{a priori} whether such a parametrization should
exist \cite{Zag09}.

In the arithmetic instances of modular origin, the sequence $c_n$ necessarily satisfies a linear differential equation \cite{Zag08}.
On the other hand, for any sequence produced by self-replication, we can (try to) guess a linear differential equation
for its generating function by computing a few initial terms.
Once guessed, we can then use closure properties of holonomic functions to prove that the
solution to the differential equation does indeed satisfy the self-replicating functional equation.
If the generating function satisfies a linear differential equation then this somewhat brute
experimental approach is guaranteed to find this equation and prove its correctness in finite time.
Showing that the sequence is not holonomic is more painful, as the time required by this brute force approach is infinite.

On the other hand, all known arithmetic instances lead to sequences which
satisfy strong divisibility properties. For instance, for primes $p$, these
sequences $c(n)$ satisfy the \emph{Lucas congruences}
\begin{equation}
c (n) \equiv c (n_0) c (n_1) \dotsb c (n_r) \pmod{p},
\label{lucas}
\end{equation}
where $n = n_0 + n_1 p + \cdots + n_r p^r$ is the expansion of $n$ in base~$p$.
In the case of Ap\'ery-like numbers, these congruences were shown in~\cite{MS16}
and further general results were recently obtained in~\cite{ABD16}.
In addition, these sequences satisfy, or are conjectured to satisfy,
the $p^{\ell r}$-congruences
\begin{equation}
c (m p^r) \equiv c (m p^{r - 1}) \pmod{p^{\ell r}}
\label{pr}
\end{equation}
with $\ell \in \{ 1, 2, 3 \}$, which in interesting cases hold for $\ell > 1$
and are referred to as \emph{supercongruences}, e.g., see \cite{OSS16} and the references therein.

In the following examples, we therefore consider self-replicating functional
equations of basic shapes, and search for values of the involved parameters
such that the corresponding sequences satisfy such congruences. In each case,
in which we observe the Lucas congruences and $p^r$-congruences, we are able to
identify these sequences and relate them to modular forms.

\medskip
The recursion \eqref{rec} together with the original self-replicating equation \eqref{alg}
suggest that the choice $\lambda^2 = \mu$ in~\eqref{alg0} is of
special arithmetic significance. This is substantiated in the next example.

\begin{example}
\label{ex1}
We consider the functional relation \eqref{alg0} with integer parameters $\lambda$
and $\mu$. The terms $c_n (\lambda, \mu)$ of the corresponding sequence are
polynomials in $\lambda$ and $\mu$. The first few are
\begin{equation*}
1, \quad 2 (\mu - \lambda), \quad 7 \mu^3 - 10 \lambda \mu + 3 \lambda^2 + 2 \mu - 2 \lambda, \quad \dotsc.
\end{equation*}
Using these initial terms we then determine congruence conditions on the
parameters $\lambda, \mu$, which need to be satisfied in order that the
$p^r$-congruences \eqref{pr}, with $\ell = 1$, hold for all primes $p$.
In the case $\lambda = \mu$ we have $c_n (\lambda, \lambda) = 0$ for $n \geq 1$,
and so we assume $\lambda \ne \mu$ in the sequel. Remarkably,
the empirical data suggests that these congruences only hold in a finite
number of cases, all of modular origin, as well as an infinite family of
unclear origin. Specifically, the $p^r$-congruences appear to hold modulo
$p^r$ for all primes $p$ only in the cases $\lambda^2 = \mu$ with
\begin{equation*}
\lambda \in \{ - 2, - 1, 2, 4, 16 \}
\end{equation*}
and the infinite family of cases $\lambda = - 2 \mu$ where $\mu$ is any even
integer. If there exist further cases, then $\min \{ | \mu |, | \lambda | \} > 10,000$
(and $| \lambda | > 100, 000$ in the particularly relevant case $\lambda^2 = \mu$).

Here is the analysis of the cases when $\lambda^2 = \mu$:
\begin{itemize}
\item
The case $\lambda = 2$ is the level 7 sequence from the introduction (and Section~\ref{sec3} below).
\item
The case $\lambda = - 2$ is the level 3 sequence $c_n=c_n(-2,4)$ with
\begin{equation*}
\sum_{n = 0}^{\infty} c_n z^n = \biggl(\sum_{n = 0}^{\infty}\binom{2n}{n}\binom{3n}{n} z^n \biggr)^2
\end{equation*}
(see Example~\ref{ex4}).
\item
The case $\lambda = - 1$ produces the sequence $c_n(-1,1) = 2^n\binom{2 n}{n}$, which clearly satisfies $p^r$-congruences (but does not
satisfy $p^{2 r}$-congruences). The corresponding generating function $\sum_{n = 0}^{\infty} c_n z^n$ is algebraic.
\item
The case $\lambda = 4$ produces the sequence $c_n=c_n(4,16)$ with
\begin{equation*}
\sum_{n = 0}^{\infty} c_n z^n = \biggl(\sum_{n = 0}^{\infty}\binom{2n}{n}^2\binom{3n}{n} z^n \biggr)^2,
\end{equation*}
which appears to satisfy $p^{3 r}$-supercongruences. The corresponding modular form has level 3 and weight~4.
\item
The case $\lambda = 16$ produces the sequence $c_n=c_n(16,256)$ with (note the exponent)
\begin{equation*}
\sum_{n = 0}^{\infty} c_n z^n = \biggl(\sum_{n = 0}^{\infty} \binom{2n}{n}\binom{3n}{n}\binom{6n}{3n} z^n \biggr)^4,
\end{equation*}
which appears to satisfy $p^{2 r}$-supercongruences. The modular form has level~1 and weight~8.
\end{itemize}

Apparently, the sequences corresponding to
$(\lambda,\mu) = (-4\alpha,2\alpha)$, that is,
$$
\frac1{(1+2\alpha z)^2}\,f\biggl(\frac z{(1+2\alpha z)^3}\biggr)
=\frac1{(1-4\alpha z)^2}\,f\biggl(\frac{z^2}{(1-4\alpha z)^3}\biggr),
$$
satisfy $p^r$ congruences for all integers $\alpha$ but do not appear to be holonomic.
However, these sequences fail to satisfy Lucas congruences (except modulo $2$ and $3$, for trivial reasons).
In the case $\alpha=1$ the first few terms of the sequence $c_n=c_n(-4,2)$ are
\begin{equation*}
1, \; 12, \; 168, \; 2496, \; 38328, \; 600672, \; 9539808, \; 152891520, \; \dotsc.
\end{equation*}
The first 250 terms do not reveal a linear recurrence with polynomial coefficients, suggesting that the sequence is not holonomic,
and we were unable to otherwise identify the sequence.  It does take the slightly simplified form
\begin{equation*}
c_n = \sum_{k = 0}^n d_k d_{n - k},
\qquad\text{where}\quad
d_n = 1, \; 6, \; 66, \; 852, \; 11874, \; 172860, \; 2586108, \; \dotsc.
\end{equation*}
Similar comments apply to the case $\alpha=-1$, which results in an alternating sequence.
\qed
\end{example}

There are many natural shapes of self-replicating functional equations that one could similarly analyze.  Here, we restrict ourselves to a variation of~\eqref{alg0}.

\begin{example}
\label{ex3}
In the spirit of the previous example, consider the functional relations
$$
\frac1{1+\mu z}\,f\biggl(\frac z{(1+\mu z)^2}\biggr)
=\frac1{1+\lambda z}\,f\biggl(\frac{z^2}{(1+\lambda z)^2}\biggr).
$$
We find supercongruences modulo $p^{3r}$, that is \eqref{pr} with $\ell=3$, in the cases $(\lambda,\mu) = (- 2, 0)$, $(\lambda,\mu) = (0, - 2)$ and $(\lambda,\mu) = (0, 4)$,
which correspond to the sequences $\binom{2 n}{n}$, $(- 1)^n\binom{2 n}{n}$ and ${\binom{2 n}{n}}^2$, respectively.

The case $(\lambda,\mu) = (-8, 16)$ produces the sequence $a_n$ with
\begin{equation*}
\sum_{n = 0}^{\infty} a_n z^n = \biggl(\sum_{n = 0}^{\infty} \binom{2n}{n}\binom{4n}{2n} z^n \biggr)^2,
\end{equation*}
which appears to satisfy $p^{2 r}$-supercongruences; this is a modular form of level 2 given in Example~\ref{ex4} below.

The cases $(\lambda,\mu) = (- 2, 0)$ and $(\lambda,\mu) = (0, - 2)$ generalize to the case $(\lambda,\mu) = (- \alpha - 2,\alpha)$,
which corresponds to the sequence $(\alpha + 1)^n \binom{2 n}{n}$. For obvious reasons we have
$p^r$-congruences for all of these, as well as Lucas congruences and recurrences.

In the cases $(\lambda,\mu) = (2 \alpha,0)$, we appear to have
$p^r$-congruences but no Lucas congruences and no recurrences.

We did not find additional instances of congruences modulo $p^r$. If
there exist any, then $\min \{ | \mu |, | \lambda | \} > 10,000$.
\qed
\end{example}

The arithmetic sequences observed in this section all satisfy strong congruences.
It would be interesting to understand to which degree such congruences
characterize these arithmetic instances of the self-replication process.

\section{Self-replicating functional equations}
\label{sec3}

For levels $\ell\in\{2,3,4,5,7\}$ we use the modular functions
$$
z_\ell(\tau)=q\biggl(\prod_{j=1}^\infty\frac{1-q^{\ell j}}{1-q^j}\biggr)^{24/(\ell-1)},
\qquad\text{where}\quad q=e^{2\pi i\tau},
$$
and the corresponding weight two modular forms
$$
P_\ell(\tau)=q\frac{\d }{\d q} \log z_\ell
$$
to define implicitly the analytic functions $f_\ell(x)$ by
\begin{gather*}
P_2(\tau) = f_2\biggl(\frac{z_2(\tau)}{1+64z_2(\tau)}\biggr),
\quad
P_3(\tau) = f_3\biggl(\frac{z_3(\tau)}{1+27z_3(\tau)}\biggr),
\quad
P_4(\tau) = f_4\biggl(\frac{z_4(\tau)}{1+16z_4(\tau)}\biggr),
\\
P_5(\tau) = f_5\biggl(\frac{z_5(\tau)}{1+22z_5(\tau)+125z_5(\tau)^2}\biggr)
\quad\text{and}\quad
P_7(\tau) = f_7\biggl(\frac{z_7(\tau)}{1+13z_7(\tau)+49z_7(\tau)^2}\biggr).
\end{gather*}
Explicitly we have
\begin{align*}
f_2(x)
&=\biggl(\sum_{n=0}^\infty\binom{2n}n\binom{4n}{2n}x^n\biggr)^2
= {}_2F_1\biggl(\begin{matrix}\frac14,\,\frac34\\1\end{matrix}\biggm|64x\biggr)^2,
\\
f_3(x)
&=\biggl(\sum_{n=0}^\infty\binom{2n}n\binom{3n}nx^n\biggr)^2
= {}_2F_1\biggl(\begin{matrix}\frac13,\,\frac23\\1\end{matrix}\biggm|27x\biggr)^2,
\\
f_4(x)
&=\biggl(\sum_{n=0}^\infty{\binom{2n}n}^2x^n\biggr)^2
= {}_2F_1\biggl(\begin{matrix}\frac12,\,\frac12\\1\end{matrix}\biggm|16x\biggr)^2,
\\ \intertext{where the standard hypergeometric notation is used, and also}
f_5(x)
&=\sum_{n=0}^\infty\biggl\{\binom{2n}n\sum_{k=0}^n{\binom nk}^2\binom{n+k}k\biggr\} x^n,
\end{align*}
while binomial expressions for the coefficients of $f_7(x)$ are given in~\eqref{eq:s7}.

\begin{example}
\label{ex4}
Each of these can be characterized by a self-replicating functional equation:
\begin{align*}
\frac{1}{1+16z}\, f_2\biggl(\frac z{(1+16z)^2}\biggr) &= \frac{1}{1-8z}\, f_2\biggl(\frac{z^2}{(1-8z)^2}\biggr),
\\
\frac{1}{(1+4z)^2}\, f_3\biggl(\frac z{(1+4z)^3}\biggr) &= \frac{1}{(1-2z)^2}\, f_3\biggl(\frac{z^2}{(1-2z)^3}\biggr),
\\
\frac{1}{(1+4z)^2}\, f_4\biggl(\frac z{(1+4z)^2}\biggr) &=  f_4(z^2),
\\
\frac{1}{1+8z}\, f_5\biggl(\frac z{(1+4z)(1+8z)^2}\biggr) &=  \frac{1}{1+2z}\, f_5\biggl(\frac{z^2}{(1+4z)(1+2z)^2}\biggr),
\\
\frac{1}{(1+4z)^2}\, f_7\biggl(\frac z{(1+4z)^3}\biggr) &= \frac{1}{(1+2z)^2}\, f_7\biggl(\frac{z^2}{(1+2z)^3}\biggr).
\end{align*}
The equations for $f_3$ and $f_7$ are discussed in Example~\ref{ex1}, and the one for $f_2$ shows up in Example~\ref{ex3}.
These instances demonstrate that the example \eqref{alg} for level 7 may be viewed as part of a bigger family.

We expect similar self-replicating equations for levels 9, 13 and 25; see \cite[Section~8]{CY15}.
\qed
\end{example}

In order to give some more weight 2 examples, define
\begin{align*}
\hat f_2(x)
&=\sum_{n=0}^\infty{\binom{2n}n}^2\binom{4n}{2n}x^n
= {}_3F_2\biggl(\begin{matrix}\frac14,\,\frac12,\,\frac34\\1,\,1\end{matrix}\biggm|256x\biggr),
\\
\hat f_3(x)
&=\sum_{n=0}^\infty{\binom{2n}n}^2\binom{3n}nx^n
= {}_3F_2\biggl(\begin{matrix}\frac13,\,\frac12,\,\frac23\\1,\,1\end{matrix}\biggm|108x\biggr),
\\
\hat f_4(x)
&=\sum_{n=0}^\infty{\binom{2n}n}^3x^n
= {}_3F_2\biggl(\begin{matrix}\frac12,\,\frac12,\,\frac12\\1,\,1\end{matrix}\biggm|64x\biggr)
\intertext{and}
\hat f_5(x)
&= \sum_{n=0}^\infty \biggl\{\sum_{k=0}^n (-1)^{n-k}{\binom nk}^3\binom{4n-5k}{3n} \biggr\}x^n.
\end{align*}

\begin{example}
\label{ex5}
A functional equation for $\hat f_2$ is given by
$$
\frac{1}{1+27z}\,\hat f_2\biggl(\frac z{(1+27z)^4}\biggr)
= \frac{1}{1+3z}\,\hat f_2\biggl(\frac{z^3}{(1+3z)^4}\biggr).
$$
This may be regarded as a cubic version of the functional equations studied in Examples~\ref{ex1}
and~\ref{ex3}.

A functional equation for $\hat f_3$ is given by the case $\lambda=4$ in Example~\ref{ex1}.

Functional equations for $\hat f_4$ and $\hat f_5$ are of a different type:
\begin{align*}
\frac{1}{1+8z}\,\hat f_4\biggl(z\biggl(\frac{1-z}{1+8z}\biggr)^3\biggr)
&=\hat f_4\biggl(z^3\biggl(\frac{1-z}{1+8z}\biggr)\biggr),
\\
\frac{1}{(1+4z)^2}\,\hat f_4\biggl(z\biggl(\frac{1-z}{1+4z}\biggr)^5\biggr)
&=\hat f_4\biggl(z^5\biggl(\frac{1-z}{1+4z}\biggr)\biggr),
\\
\frac{1}{1-5z}\,\hat f_5\biggl(z\biggl(\frac{1-z}{1-5z}\biggr)^2\biggr)
&=\hat f_5\biggl(z^2\biggl(\frac{1-z}{1-5z}\biggr)\biggr).
\qed
\end{align*}
\end{example}

Example~\ref{ex5} suggests considering functional equations of the type
$$
t(z)\, f\bigl(r(z)\,s^n(z)\bigr) =  f\bigl(r^n(z)\,s(z)\bigr)
$$
where $r(z)$, $s(z)$ and $t(z)$ are rational functions with
$$
r(z) = z+O(z^2), \quad s(z)=1+O(z), \quad t(z) = 1+O(z),
$$
and $n\geq 2$ is an integer.

\begin{example}
\label{ex7}
We list three more instances that may be compared with
the functional equation for $f_4$ in Example~\ref{ex4}:
\begin{align*}
g_b(z^2)&=\frac{1}{1+3z}\,g_b\biggl(z\,\frac{1-z}{1+3z}\biggr),
\\
g_c(z^3)&=\frac{1}{1+2z+4z^2}\,g_c\biggl(z\,\frac{1-z+z^2}{1+2z+4z^2}\biggr)
\\ \intertext{and}
g_5(z^5)&=\frac{1}{1+3z+4z^2+2z^3+z^4}\,g_5\biggl(z\,\frac{1-2z+4z^2-3z^3+z^4}{1+3z+4z^2+2z^3+z^4}\biggr).
\end{align*}
The solutions to these functional equations that satisfy the initial condition $g(0)=\nobreak1$ are given by
$$
g_b(x) = \sum_{n=0}^\infty \biggl\{\sum_{k=0}^\infty {\binom{n}{k}}^2 \binom{2k}{k}\biggr\}x^n,
\qquad
g_c(x) = \sum_{n=0}^\infty \biggl\{\sum_{k=0}^\infty {\binom{n}{k}}^3 \biggr\}x^n
$$
and
$$
g_5(x) = \sum_{n=0}^\infty \biggl\{\sum_{k=0}^\infty {\binom{n}{k}}^2 \binom{n+k}{k}\biggr\}x^n.
$$
The functions $g_b$ and $g_c$ may be parameterized by level~6 modular forms, while~$g_5$ may be parameterized
by level~5 modular forms, e.g., see \cite[Theorem~3.4 and Tables~1 and~2]{CC12}.
The functions $g_b$, $g_c$ and $g_5$ are denoted by (c), (a) and (b), respectively, in \cite[Eq.~(28)]{AvSZ11}.
They are generating functions for the sequences C, A and D, respectively, in~\cite{Zag09}.
\qed
\end{example}

\section({Formulas for 1/\003\300 and AGM-type iterations}){Formulas for $1/\pi$ and AGM-type iterations}
\label{sec4}

\begin{example}
\label{ex8}
We can use the transformation \eqref{alg} as $z\to(1/8)^-$, so that
$z/(1+4z)^3\to(1/27)^-$ and $z/(1+4z)^3\to1/125$.
Applying the asymptotics and techniques from \cite[Example 5]{GZ13}
we obtain the formula
\begin{equation}
\label{N}
\sum_{n=0}^\infty u_n(4+21n)\times\frac1{5^{3n+3}}=\frac{1}{8\pi}.
\end{equation}
This can also be obtained using modular forms, e.g., see \cite[Eq.~(37) and Table 1]{Coo12}.

Consider the limiting case $h\to h_0^+$, where $h_0=-(5+\sqrt{21})/2$, of the transformation
\cite[Lemmas 4.1 and 4.3]{CC08}
\begin{align*}
&
\frac1{\sqrt{1+13h+49h^2}}\sum_{n=0}^\infty u_n\biggl(\frac h{1+13h+49h^2}\biggr)^n
\\ &\quad
=\frac1{\sqrt{1+5h+h^2}}
\,{}_3F_2\biggl(\begin{matrix} \frac16, \, \frac12, \, \frac56 \\ 1, \, 1 \end{matrix}\biggm|
\frac{1728h^7}{(1+13h+49h^2)(1+5h+h^2)^3}\biggr),
\end{align*}
so that the argument of the $_3F_2$ series on the right-hand side tends to $-\infty$,
and follow \cite[Theorem~2 and Example~2]{GZ13} to obtain
\begin{equation}
\label{N21a}
\sum_{n=0}^\infty u_n\biggl(\frac{3\sqrt{21}-14}{56}\biggr)^n
=\biggl(\frac{128(\sqrt{7}-\sqrt{3})}{49(5-\sqrt{21})^2}\biggr)^{1/3}\frac{\sqrt\pi}{3\Gamma(\frac56)^3}
\end{equation}
and
\begin{equation}
\label{N21}
\sum_{n=0}^\infty u_n\bigl(6\sqrt{21}-20+15(\sqrt{21}-2)n\bigr)\biggl(\frac{3\sqrt{21}-14}{56}\biggr)^n=\frac{8\sqrt7}\pi.
\end{equation}
The series \eqref{N21} corresponds to the case $N=21$ of
the identity (39) in \cite{Coo12}.
\qed
\end{example}

A principal feature of an AGM-type algorithm is a sequence of algebraic approximations
$a_k$, where $k=0,1,2,\dots$, that converge to a given number $\xi$
at a certain rate $r>1$, so that $|\xi-a_{k+1}|<C|\xi-a_k|^r$ for some absolute constant $C>0$.
If $r=2$, as for the classical AGM algorithm of Brent~\cite{Bre76} and Salamin~\cite{Sal76}, then it is said to be a
quadratic iteration. The next example uses a functional equation to
produce quadratic AGM-type iterations to numbers $\xi$ of the form
\begin{equation}
\sum_{n=0}^\infty u_n\,(a_0+b_0n)x_0^n=\xi,
\label{id0}
\end{equation}
where $a_0$, $b_0$ and $x_0$ are fixed real algebraic numbers, where the sequence $u_n$ is given in~\eqref{eq:s7}.
The iterations are based on the self-replicating equation \eqref{alg}.

\begin{example}
\label{ex9}
Let $a_0$, $b_0$ and $x_0$ be chosen in accordance with \eqref{id0}, while
$z_0$ is the real solution to
$$
\frac{z}{(1+4z)^3}=x_0
$$
of the smallest absolute value.
Define the sequences $a_k$, $b_k$ and $z_k$ for all $k\ge1$ recursively as follows:
\begin{align*}
a_{k+1}&=a_k\,\frac{(1+4z_k)^2}{(1+2z_k)^2}+b_k\,\frac{4z_k(1+4z_k)^2}{(1+2z_k)^3(1-8z_k)},
\\
b_{k+1}&=2b_k\,\frac{(1+4z_k)^3(1-z_k)}{(1+2z_k)^3(1-8z_k)},
\\
z_{k+1}&=\frac{2z_k^2}{1+6z_k+(1+2z_k)\sqrt{1+8z_k}}.
\end{align*}
Then $a_k$ converges to $\xi$ in~\eqref{id0} quadratically.

\medskip
To justify the above iteration, we couple the self-replicating equation \eqref{alg}
with the one obtained from it by differentiation:
\begin{multline}
\sum_{n=0}^\infty u_n\biggl(\biggl(A-\frac{8Bz}{1+4z}\biggr)+nB\biggl(1-\frac{12z}{1+4z}\biggr)\biggr)\frac{z^k}{(1+4z)^{3k+2}}
\\
=\sum_{n=0}^\infty u_n\biggl(\biggl(A-\frac{4Bz}{1+2z}\biggr)+nB\biggl(2-\frac{6z}{1+2z}\biggr)\biggr)\frac{z^{2k}}{(1+2z)^{3k+2}}.
\label{alg1}
\end{multline}
In addition to the above three sequences we also consider
\begin{equation}
\label{x-r}
x_k=\frac{z_k}{(1+4z_k)^3};
\end{equation}
note that the definition of $z_{k+1}$ implies that
$$
x_{k+1}=\frac{z_{k+1}}{(1+4z_{k+1})^3}
=\frac{z_k^2}{(1+2z_k)^3}.
$$
Then induction on $k$ then shows that
\begin{equation}
\sum_{n=0}^\infty u_n\,(a_k+b_kn)x_k^n=\xi.
\label{idk}
\end{equation}
Indeed, on the $k$-th step we first define the numbers $A$, $B$ and $z=z_k$ such that
$$
\frac{z}{(1+4z)^3}=x_k,
\quad
\frac{B}{(1+4z)^2}\biggl(1-\frac{12z}{1+4z}\biggr)=b_k,
\quad
\frac1{(1+4z)^2}\biggl(A-\frac{8Bz}{1+4z}\biggr)=a_k
$$
and assign to the data
\begin{equation}
\label{x-iteration}
x_{k+1}=\frac{z^2}{(1+2z)^3},
\quad
b_{k+1}=\frac{B}{(1+2z)^2}\biggl(2-\frac{6z}{1+2z}\biggr)
\end{equation}
and
\begin{equation}
\label{a-iteration}
a_{k+1}=\frac1{(1+2z)^2}\biggl(A-\frac{4Bz}{1+2z}\biggr).
\end{equation}
Application of \eqref{alg1} implies the validity of \eqref{idk} for $k$ replaced with $k+1$;
taking $z=z_k$ and eliminating $A,B$ lead to the formulas for $a_{k+1}$, $b_{k+1}$ and $x_{k+1}$
by means of $a_k$, $b_k$ and $z_k$ only.

In this construction for $x_k$ sufficiently close to 0 we have $x_{k+1}\approx z_k^2\approx x_k^2$,
so that $x_k$ tends to~0 quadratically. Performing the limit in \eqref{idk} as $k\to\infty$ we obtain
$$
\xi=\lim_{k\to\infty}\sum_{n=0}^\infty u_n\,(a_k+b_kn)x_k^n=\lim_{k\to\infty}a_k.
$$
Furthermore, $|\xi-a_k|$ is roughly of magnitude $x_k$, so that $a_k$ converges to $\xi$ quadratically.

\medskip
To execute the algorithm one can use \eqref{N}
to produce the initial conditions
\begin{equation}
\label{i-c}
a_0=\frac{4}{125},\quad b_0=\frac{21}{125},\quad x_0=\frac{1}{125}, \quad z_0=\biggl(\frac{\sqrt2-1}2\biggr)^3=0.00888347\dotsc.
\end{equation}
Then for the iteration scheme given by \eqref{x-iteration} and \eqref{a-iteration} we have
$$
\lim_{n\to\infty} a_n = \frac{1}{2\pi}
$$
and the convergence is quadratic.

Moreover, the formula \eqref{N21a} leads to the initial conditions
$$
a_0=3\biggl(\frac{128(\sqrt{7}-\sqrt{3})}{49(5-\sqrt{21})^2}\biggr)^{-1/3}, \quad b_0=0,
\quad x_0=\frac{3\sqrt{21}-14}{56},
$$
while $z_0$ is a real solution of the equation $x_0(1+4z)^3-z=0$,
for which the iteration scheme given by \eqref{x-iteration} and \eqref{a-iteration}
has the property that
$$
\lim_{n\to\infty} a_n = \frac{\pi^{1/2}}{\Gamma(5/6)^3}=\frac{\Gamma(1/6)^3}{8\pi^{5/2}}
$$
with quadratic convergence. In a similar way, the identity \eqref{N21} leads to a different set of initial
conditions such the the iteration scheme given by \eqref{x-iteration} and \eqref{a-iteration} generates
another sequence that converges to $1/(2\pi)$.
\qed
\end{example}

Our next example makes use of the identities from Example~\ref{ex5} to produce an iteration scheme
that converges quintically to $1/\pi$.

\begin{example}
\label{ex10}
Define sequences $a_k$, $b_k$, $x_k$ and $z_k$ as follows.
Let
$$
a_0=\frac14,\quad b_0=1,\quad x_0=-\frac1{64}.
$$
Recursively, let $z_k$ be the real solution to
$$
z^{1/5}\frac{1-z}{1+4z}=x_k^{1/5}
$$
that is closest to the origin,
and then let
\begin{align*}
a_{k+1}&=a_k\, (1+4z_k)^2 + 8\,b_k\,\frac{z_k(1-z_k)(1+4z_k)^2}{1-22z_k-4z_k^2},
\\
b_{k+1}&=5\,b_k\,(1+4z_k)^2\frac{1+2z_k-4z_k^2}{1-22z_k-4z_k^2}
\end{align*}
and
$$
x_{k+1}=z_k^5\frac{1-z_k}{1+4z_k}.
$$
Then
$$
\lim_{k\to\infty} a_k = \frac{1}{2\pi}
$$
and the rate of convergence is of order 5.

\medskip
To justify the above algorithm, we differentiate the last formula of Example~\ref{ex5} to get
\begin{multline*}
\sum_{n=0}^\infty {\binom{2n}n}^3 \biggl(\biggl(A-\frac{8Bx}{1+4x}\biggr)+Bn\,\frac{1-22x-4x^2}{(1-x)(1+4x)}\biggr)
\frac{x^n(1-5x)^{5n}}{(1+4x)^{5n+2}}
\\
= \sum_{n=0}^\infty {\binom{2n}n}^3 \biggl(A+5Bn\,\frac{1+2x-4x^2}{(1-x)(1+4x)}\biggr) \frac{x^{5n}(1-x)^n}{(1+4x)^n}.
\end{multline*}
Consider an iteration scheme that goes as follows.
Suppose $x_k$, $a_k$ and $b_k$ are given, for some $k$. Compute $z$, $B$ and $A$ from the formulas
$$
z\biggl(\frac{1-z}{1+4z}\biggr)^5=x_k,
\quad
B\,\frac{1-22x-4x^2}{(1+4x)^3(1-x)}=b_k
$$
and
$$
\frac{1}{(1+4x)^2}\biggl(A-\frac{8Bx}{1+4x}\biggr)=a_k.
$$
Then, let $x_{k+1}$, $b_{k+1}$ and $a_{k+1}$ be defined by
$$
x_{k+1}=\frac{x^5(1-x)}{(1+4x)},\quad b_{k+1}=\frac{5B(1+2z-4z^2)}{(1-z)(1+4z)},\quad a_{k+1}=A.
$$
These iterations are equivalent to the ones stated in the theorem.

To get initial conditions, consider Bauer's series \cite[Table 6, $N=2$]{CC12}, \cite{Zud08}
$$
\sum_{n=0}^\infty {\binom{2n}n}^3\biggl(\frac{1}{4}+n\biggr)\biggl(-\frac1{64}\biggr)^n = \frac{1}{2\pi}.
$$
This gives $a_0=1/4$, $b_0=1$ and $x_0=-1/64$, as well as the limit
$$
\lim_{n\to \infty} a_n = \frac{1}{2\pi}.
$$
The data
\begin{alignat*}{2}
\biggl| a_0 - \frac{1}{2\pi}\biggr| &\approx 9.08 \times 10^{-2}, &\quad
\biggl| a_1 - \frac{1}{2\pi}\biggr| &\approx 6.65 \times 10^{-9}, \\
\biggl| a_2 - \frac{1}{2\pi}\biggr| &\approx 8.25 \times 10^{-47}, &\quad
\biggl| a_3 - \frac{1}{2\pi}\biggr| &\approx 4.57 \times 10^{-239}
\end{alignat*}
is good supporting evidence.

One can also use the series \cite[Table 6, $N=3$]{CC12}
$$
\sum_{n=0}^\infty {\binom{2n}n}^3\biggl(\frac16+n\biggr)\biggl(\frac{1}{256}\biggr)^n = \frac{2}{3\pi}
$$
to produce initial conditions $a_0=1/6$, $b_0=1$, $x_0=1/256$ and limiting value
$$
\lim_{n\to \infty} a_n = \frac{2}{3\pi}.
$$
Another option is to use the series \cite[Table 6, $N=7$]{CC12}
$$
\sum_{n=0}^\infty {\binom{2n}n}^3\biggl(\frac{5}{42}+n\biggr)\biggl(\frac{1}{4096}\biggr)^n = \frac{8}{21\pi},
$$
to give $a_0=5/42$, $b_0=1$, $x_0=1/4096$ and the limit
$$
\lim_{n\to \infty} a_n = \frac{8}{21\pi}.
$$
These have all been tested numerically. There is no practical advantage to starting with different initial conditions,
as the speed of convergence is determined by the quintic nature of the algorithm.
\qed
\end{example}

Observe that the sequence $b_k$ in Example~\ref{ex9} with initial conditions given by~\eqref{i-c} has the property
$$
\frac{b_k}{b_0} = 2^k  \times \sqrt{\frac{1-26x_k-27x_k^2}{1-26x_0-27x_0^2}}
$$
where $x_k$ is given by \eqref{x-r}. And, for the sequences $b_k$ in Example~\ref{ex10} we have
$$
\frac{b_k}{b_0} = 5^k\times \sqrt{\frac{1-64x_k}{1-64x_0}}.
$$
The factors $2^k$ and $5^k$ reflect the quadratic and quintic rates of convergence of the respective iterations,
and the factors $1-26x-27x^2$ and $1-64x$ are the leading coefficients in the respective differential equations
for the functions~$f(x)$, e.g., see \cite[p. 176]{Coo12}.

\section{Self-replication of Legendre polynomials}
\label{sec5}

Several generating functions of the Legendre polynomials
$$
P_n(x)={}_2F_1\biggl(\begin{matrix} -n, \, n+1 \\ 1 \end{matrix}\biggm| \frac{1-x}2 \biggr)
=\frac1{2^n}\sum_{k=0}^n {\binom nk}^2(x-1)^k(x+1)^{n-k}
$$
are known but the one we are going to treat in this section is given by
$$
F(x;z)=\sum_{n=0}^\infty{\binom{2n}n}^2P_n(x)z^n.
$$
The Bailey--Brafman identity \cite{CWZ13} allows us to recast the two-variable function as a product of two hypergeometric functions,
$$
F\biggl(\frac{U+V-2UV}{U-V};\frac{U-V}{16}\biggr)
= {}_2F_1\biggl(\begin{matrix}\frac12,\,\frac12\\1\end{matrix}\biggm|U\biggr)
\cdot {}_2F_1\biggl(\begin{matrix}\frac12,\,\frac12\\1\end{matrix}\biggm|V\biggr),
$$
and with the help of the quadratic transformation
$$
{}_2F_1\biggl(\begin{matrix}\frac12,\,\frac12\\1\end{matrix}\biggm|u^2\biggr)
=\frac1{1+u}\cdot{}_2F_1\biggl(\begin{matrix}\frac12,\,\frac12\\1\end{matrix}\biggm|\frac{4u}{(1+u)^2}\biggr)
$$
(essentially, the equation for $f_4$ in Example~\ref{ex4})
we get the following self-replicating identity:
\begin{multline}
F\biggl(\frac{u^2+v^2-2u^2v^2}{u^2-v^2};\frac{u^2-v^2}{16}\biggr)
\\
=\frac1{(1+u)(1+v)}\,F\biggl(\frac{(1+uv)(u+v)-4uv}{(1-uv)(u-v)};\frac{(1-uv)(u-v)}{4(1+u)^2(1+v)^2}\biggr).
\label{leg}
\end{multline}
Because of the way the two parameters $u$ and $v$ are tangled in the relation \eqref{leg}, it is not easy to perform
an analysis of the equation analogous to the one we have had for~\eqref{alg}. Notice, however, that $F(x;z)$ is
a particular case of Appell's $F_4$ function and that there is a similar-looking identity \cite{KS07} for Appell's $F_1$ function
(see also \cite[Chap.~4]{Sch14}). As the latter identity is known to be linked to a three-term generalization
of the arithmetic-geometric mean, we believe that \eqref{leg} and the like equations may find some interesting
number-theoretic applications.

\end{document}